\documentclass[12pt]{article}
\usepackage{amssymb, amsmath, amsfonts,dsfont, mathrsfs,euscript,eufrak,colonequals}
\usepackage{textcomp}
\usepackage{latexsym}

\newtheorem{Theorem}{Theorem}

\newcommand{\newparagraph}{\hspace{\parindent}}

\newcommand{\N}{\mathbb N}
\newcommand{\R}{\mathbb R}

\newcommand{\Expec}{ \mathds{E}\,}

\newcommand{\dd}{\,\mathrm{d}}
\newcommand{\id}{\mathds{1}}

\newcommand{\e}{\varepsilon}
\newcommand{\AppClass}{\mathcal{A}}

\newcommand{\nlambda}{\bar{\lambda}}

\newcommand{\CorrFunc}{\mathcal{K}}

\newcommand{\SEElem}{\mathbb{G}}
\newcommand{\KorElem}{\mathbb{B}}
\begin{document}
\author{A. A. Khartov}
\title{Simplified criterion of quasi-polynomial tractability\\ and its applications\footnote{The work was supported by the RFBR grant 13-01-00172, by the SPbGU grant 6.38.672.2013, and by the grant of Scientific school NSh-2504.2014.1.}}

\maketitle
\begin{abstract}
We study approximation properties of sequences of centered random elements $X_d$, $d\in\N$, with values in separable Hilbert spaces. We focus on sequences of tensor product-type random elements, which have covariance operators of corresponding tensor product form. The average case approximation complexity $n^{X_d}(\e)$ is defined as the  minimal number of evaluations of arbitrary linear functionals that is needed to approximate $X_d$ with relative $2$-average error not exceeding a given threshold $\e\in(0,1)$. The growth of $n^{X_d}(\e)$ as a function of $\e^{-1}$ and $d$ determines whether a sequence of corresponding approximation problems for $X_d$, $d\in\N$, is tractable or not. Different types of tractability were studied in the paper by M. A. Lifshits, A. Papageorgiou and H. Wo\'zniakowski (2012), where for each type the necessary and sufficient conditions were found in terms of the eigenvalues of the marginal covariance operators. We revise the criterion of quasi-polynomial tractability and provide its simplified version. We illustrate our result by applying to random elements corresponding to tensor products of squared exponential kernels. Also we extend recent result of G. Xu (2014) concerning weighted Korobov kernels.
\end{abstract}

\section{Introduction}\newparagraph
Let $X_d$, $d\in\N$, be a sequence of random elements of some normed spaces $(Q_d, \|\cdot\|_{Q_d})$, $d\in\N$, respectively, where every $Q_d$ be a space of functions of $d$ variables. How do approximation properties of $X_d$ depend on $d$? More formally, we consider the following \textit{linear tensor product approximation problems} in \textit{average case setting} (\textit{approximation problems} for short, see \cite{NovWoz1}--\cite{NovWoz3} and \cite{Rit}). Suppose that $Q_d=\otimes_{j=1}^d Q_{1,j}$ in appropriate sense, where $Q_{1,j}$, $j\in\N$, are some normed spaces. Suppose that every $Q_d$-valued random element $X_d$ is centered and has the covariance operator $K^{X_d}$ of the appropriate tensor product form $K^{X_d}=\otimes_{j=1}^d K^{X_{1,j}}$, $d\in\N$, where $K^{X_{1,j}}$ is the covariance operator of a given $Q_{1,j}$-valued centered random element $X_{1,j}$, $j\in\N$. Such  $X_d$ is called the \textit{tensor product of} $X_{1,1},\ldots, X_{1,d}$.  We approximate every $X_d$ by the finite rank sums $\widetilde{X}_d^{(n)}=\sum_{k=1}^{n} l_k (X_d) \psi_k$, where $\psi_k$ are deterministic elements of $Q_d$ and $l_k$ are continuous linear functionals from the dual space $Q_d^*$. We consider the \textit{average case approximation complexity} $n^{X_d}(\e)$ as a characteristic of approximation of the random element $X_d$. It is defined as the minimal suitable value of $n$ needed to make the relative average approximation error $\bigl(\Expec\|X_d-\widetilde{X}_d^{(n)}\|_{Q_d}^2/\Expec \|X_d\|_{Q_d}^2\bigr)^{1/2}$ smaller than a given error threshold $\e$ by choosing optimal $\psi_k$ and $l_k$ (see \cite{TraubWasWoz}). 

It is important to study tractability of the described multivariate approximation problems. Namely, the approximation complexity $n^{X_d}(\e)$ is considered as a function of two variables $\e\in(0,1)$ and $d\in\N$.  A sequence of approximation problems for $X_d, d\in\N$ is called \textit{weakly tractable} if  $n^{X_d}(\e)$ is not exponential in $d$ or/and $\e^{-1}$. Otherwise, the sequence of the problems is \textit{intractable}. Special subclasses of weakly tractable problems are distinguished depending on the types of majorants for the quantity $n^{X_d}(\e)$ \textit{for all} $d\in\N$ and $\e\in(0,1)$. For example, the sequence of approximation problems for $X_d$, $d\in\N$, is called \textit{polynomially tractable} if the majorant of $n^{X_d}(\e)$ is of order $\e^{-s}d^p$ with some non-negative constants $s$ and $p$. In the case $p=0$ the sequence of the problems  is called \textit{strongly polynomially tractable}. \textit{Quasi-polynomial tractability}, which was introduced in \cite{GnewWoz}, means that a majorant of $n^{X_d}(\e)$ is of order $\exp\{s(1+|\ln \e|)(1+\ln d)\}$ with some constant $s\geqslant 0$. In the recent paper \cite{LifPapWoz1} these types of  tractability of the described approximation problems were investigated for separable Hilbert spaces $Q_{1,j}$, $j\in\N$. For each tractability type the necessary and sufficient conditions were found in terms  of eigenvalues  of the marginal covariance operators $K^{X_{1,j}}$, $j\in\N$ (the asymptotic setting ``$\e$ is fixed, $d\to\infty$'' was considered in \cite{Khart}, \cite{LifTul}, and \cite{LifZani}). 

However it is seems that criterion of quasi-polynomial tractability from \cite{LifPapWoz1} was formulated in unfinished form. It is hardly applicable to concrete examples of $(X_d)_{d\in\N}$, because it usually requires a lot of additional calculations. The aim of this paper is to provide  a simplified criterion of quasi-polynomial tractability, which will be more convenient for applications.

The paper is organized as follows. In Section 2 we provide necessary definitions and facts  concerning  linear tensor product approximation problems defined over Hilbert spaces. In Section 3 we obtain a new criterion of quasi-polynomial tractability for these problems. In Section 4 for illustration we apply this criterion to well known example. Namely, we consider random elements corresponding to tensor products of weighted Korobov kernels with varying positive weight parameters $g_j\leqslant 1$ and smoothness parameters $r_j>1/2$, $j\in\N$. We show that assumptions on monotonicity of the sequences $(g_j)_{j\in\N}$ and $(r_j)_{j\in\N}$ can be omitted. Thus we extend the corresponding results from \cite{LifPapWoz1} and \cite{Xu}. We also apply our criterion to random elements corresponding to tensor products of squared exponential kernels with varying length scales.

Throughout the article, we use the following notation. We write $a_n\asymp b_n$ whenever there are positive constants $c_1$ and $c_2$ such that $c_1 b_n\leqslant a_n\leqslant c_2 b_n$ for all $n$. We denote by $\N$ and $\R$ the sets of positive integers and real numbers, respectively. We set $\ln_+ x\colonequals \max\{1, \ln x\}$ for all $x>0$. The quantity $\id(A)$ equals one for the true logic propositions $A$ and zero for the false ones. We always use $\|\cdot\|_B$ for the norm, which some space $B$ is equipped with.

\section{Linear tensor product approximation problems}\newparagraph
Suppose that we have a sequence of zero-mean random elements $X_{1,j}$, $j\in\N$, of separable Hilbert spaces $H_{1,j}$, $j\in\N$, respectively. We always assume that every $X_{1,j}$ satisfies $\Expec \|X_{1,j}\|^2_{H_{1,j}}<\infty$. i.e. the covariance operator of $X_{1,j}$, denoted by $K^{X_{1,j}}$, has the finite trace. Consider the sequence $X_d$, $d\in\N$, of increasing \textit{tensor products} of random elements $X_{1,j}$, $j\in\N$. Namely, every $X_d$ is a zero-mean random element of the Hilbertian tensor product $H_d\colonequals \otimes_{j=1}^d H_{1,j}$, with the covariance operator $K^{X_d}\colonequals\otimes_{j=1}^d K^{X_{1,j}}$, $d\in\N$. Following \cite{KarNazNik} and \cite{LifZani}, for random element $X_d$ of such type  we use the notation $X_d=\otimes_{j=1}^d X_{1,j}$ for short.

We will investigate the \textit{average case approximation complexity} (simply the \textit{approximation complexity} for short) of $X_d$, $d\in\N$:
\begin{eqnarray}\label{def_nXde}
n^{X_d}(\e)\colonequals\min\bigl\{n\in\N:\, e^{X_d}(n)\leqslant \e\, e^{X_d}(0)  \bigr \},
\end{eqnarray}
where $\e\in(0,1)$ is a given error threshold, and
\begin{eqnarray*}
e^{X_d}(n)\colonequals\inf\Bigl\{ \bigl(\Expec\bigl\|X_d - \widetilde X^{( n)}_d\bigr\|_{H_d}^2\bigr)^{1/2} :  \widetilde X^{(n)}_d\in \AppClass_n^{X_d}\Bigr\}
\end{eqnarray*}
is the smallest 2-average error among all linear approximations of $X_d$, $d\in\N$, having rank $n\in\N$. The corresponding classes  of linear algorithms are denoted by
\begin{eqnarray*}
\AppClass_n^{X_d}\colonequals \Bigl\{\sum_{m=1}^{n}l_m(X_d)\,\psi_m :  \psi_m \in H_d,\, l_m\in H_d^*  \Bigr \}, \quad d\in\N,\quad n\in\N.
\end{eqnarray*}
We always work with \textit{relative errors}, thus taking into account the following ``size'' of $X_d$:
\begin{eqnarray*}
e^{X_d}(0)\colonequals  \bigl(\Expec\|X_d\|_{H_d}^2\bigr)^{1/2}<\infty,
\end{eqnarray*}
which is the approximation error of $X_d$ by zero element of $H_d$. 

The approximation complexity $n^{X_d}(\e)$ is considered as a function depending on two variables $d\in\N$ and $\e\in(0,1)$. According to \cite{NovWoz1}, a sequence of approximation problems for $X_d$, $d\in\N$, is  called
\begin{itemize}
\item \textit{weakly tractable} if
\begin{eqnarray}\label{def_WT}	
\lim_{d+\e^{-1}\to\infty}\dfrac{\ln n^{X_d}(\e)}{d+\e^{-1}}=0;
\end{eqnarray}
\item \textit{quasi-polynomially tractable} if there are numbers $C>0$ and $s\geqslant0$ such that
\begin{eqnarray}\label{def_QPT}	
n^{X_d}(\e)\leqslant C\exp\bigl\{s(1+\ln \e^{-1})(1+\ln d)\bigr\}\quad\text{for all}\quad d\in\N, \,\,\e\in(0,1);
\end{eqnarray}
\item \textit{polynomially tractable} if there are numbers  $C>0$, $s\geqslant 0$, and $p\geqslant0$ such that 
\begin{eqnarray}\label{def_PT}	
n^{X_d}(\e)\leqslant C\,\e^{-s}\,d^{\,p} \quad\text{for all}\quad d\in\N, \,\,\e\in(0,1);
\end{eqnarray}
\item \textit{strong polynomially tractable} if there are numbers $C>0$ and $s\geqslant0$ such that
\begin{eqnarray}\label{def_SPT}	
n^{X_d}(\e)\leqslant C\,\e^{-s}\quad\text{for all}\quad d\in\N, \,\,\e\in(0,1).
\end{eqnarray}
\end{itemize}
If the sequence of approximation problems is not weakly tractable, then it is called \textit{intractable}. 

Let $(\lambda^{X_d}_k)_{k\in\N}$ and $(\psi^{X_d}_k)_{k\in\N}$ denote the non-increasing sequence of eigenvalues and the corresponding  sequence of eigenvectors of $K^{X_d}$, respectively, i.e. $K^{X_d} \psi^{X_d}_k=\lambda^{X_d}_k\psi^{X_d}_k$, $k\in\N$. If $X_d$ is a random element of $p$-dimensional space, then we formally set $\lambda_k^{X_d}\colonequals0$, and $\psi_k^{X_d}\colonequals0$ for $k>p$. Let $\Lambda^{X_d}$ denote the trace of $K^{X_d}$, i.e.
\begin{eqnarray}\label{def_LambdaXd}
\Lambda^{X_d}\colonequals\sum_{k=1}^\infty \lambda^{X_d}_k=\Expec \|X_d\|_{H_d}^2=e^{X_d}(0)^2<\infty,\quad d\in\N.
\end{eqnarray} 

It is well known (see \cite{WasWoz}) that for any $n\in\N$ the following random element
\begin{eqnarray}\label{def_Xdn}
\widetilde X^{(n)}_d\colonequals\sum_{k=1}^n (X_d,\psi^{X_d}_k)_{H_d}\, \psi^{X_d}_k\in\AppClass_n^{X_d}
\end{eqnarray}
minimizes the 2-average case error. Hence formula \eqref{def_nXde} is reduced to
\begin{eqnarray*}
n^{X_d}(\e)=\min\Bigl\{n\in\N:\, \Expec 
\bigl\|X_d-\widetilde{X}^{(n)}_d\bigr\|_{H_d}^2\leqslant\e^2\,\Expec \|X_d\|_{H_d}^2 \Bigr\},\quad  d\in\N,\,\, \e\in(0,1).
\end{eqnarray*}
On account of \eqref{def_LambdaXd}, \eqref{def_Xdn}, and $\Expec (X_d,\psi^{X_d}_k)_{H_d}^2=\lambda^{X_d}_k$, $k\in\N$,  we infer the following representation of the approximation complexity:
\begin{eqnarray*}
n^{X_d}(\e)=\min\Bigl\{n\in\N:\, \sum_{k=n+1}^\infty \lambda^{X_d}_k\leqslant\e^2\,\Lambda^{X_d} \Bigr\},\quad  d\in\N,\,\, \e\in(0,1).
\end{eqnarray*}

Due to the tensor structure of $K^{X_d}\colonequals \otimes_{j=1}^d K^{X_{1,j}}$, $(\lambda^{X_d}_k)_{k\in\N}$ is the ordered sequence of the following products
\begin{eqnarray*} 
\prod_{j=1}^d\lambda^{X_{1,j}}_{k_j},\qquad k_1, k_2, \ldots,k_d\in\N,
\end{eqnarray*}
where every $(\lambda_k^{X_{1,j}})_{k\in\N}$ is the non-increasing sequence of eigenvalues of $K^{X_{1,j}}$, $j\in\N$. Here if some $X_{1,j}$ is a random element of $p$-dimensional space, then we formally set $\lambda^{X_{1,j}}_k=0$ for $k>p$. Of course, we always assume that $\lambda^{X_{1,j}}_1>0$ for all $j\in\N$.  Let $\Lambda^{X_{1,j}}$ be the trace of $K^{X_{1,j}}$, i.e. 
\begin{eqnarray*}
\Lambda^{X_{1,j}}=\sum_{k=1}^\infty \lambda^{X_{1,j}}_k=\Expec\|X_{1,j}\|^2_{H_{1,j}}<\infty, \quad j\in\N.
\end{eqnarray*}
Hence for every $\Lambda^{X_d}$, $d\in\N$, we have the formula
\begin{eqnarray}\label{conc_LambdaXd}
\Lambda^{X_d}=\sum_{k_1, k_2, \ldots,k_d\in\N} \prod_{j=1}^d\lambda^{X_{1,j}}_{k_j}=\prod_{j=1}^d\sum_{i=1}^\infty\lambda^{X_{1,j}}_i=\prod_{j=1}^d\Lambda^{X_{1,j}},\quad d\in\N.
\end{eqnarray}

Thus the tractability for $X_d=\otimes_{j=1}^d X_{1,j}$, $d\in\N$, can be fully characterized by the sequences $(\lambda_k^{X_{1,j}})_{k\in\N}$, $j\in\N$. It was done for described tractability types \eqref{def_WT}--\eqref{def_SPT} in the recent paper \cite{LifPapWoz1}. In the next section we focus only on the one of those types.

\section{Quasi-polynomial tractability}\newparagraph
Here we study quasi-polynomial tractability of linear tensor product approximation problems, which were described in the previous section. From now on, we use the notation $\nlambda^{X_{1,j}}_k\colonequals \lambda^{X_{1,j}}_k/\Lambda^{X_{1,j}}$, $k\in\N$,  $j\in\N$. Thus
\begin{eqnarray}\label{cond_sumnlambdaX1jk}
\sum_{k=1}^{\infty} \nlambda^{X_{1,j}}_k=1,\quad j\in\N.
\end{eqnarray}
In the next theorems  we adopt the convention $|\ln 0\,|\cdot 0=0$.

We first recall the criterion of quasi-polynomial tractability that was obtained by  M. A. Lifshits, A. Papageorgiou and H. Wo\'zniakowski in the paper \cite{LifPapWoz1}.
\begin{Theorem}\label{th_oldQPTcrit}
The sequence of approximation problems for $X_d=\otimes_{j=1}^d X_{1,j}$, $d\in\N$, is quasi-polynomially tractable iff
\begin{eqnarray}\label{th_oldQPTcrit_cond}
\sup_{d\in\N} \prod_{j=1}^{d} \sum_{k=1}^{\infty} \bigl(\nlambda^{X_{1,j}}_k\bigr)^{1-\tfrac{\gamma}{\ln_+ d}}<\infty
\end{eqnarray}
for some $\gamma\in(0,1)$. If we have quasi-polynomially tractability then
\begin{eqnarray}\label{th_oldQPTcrit_necesscond}
\sup_{d\in\N} \dfrac{1}{\ln_+ d}\sum_{j=1}^{d} \sum_{k=1}^{\infty} |\ln\nlambda^{X_{1,j}}_k|\,\nlambda^{X_{1,j}}_k<\infty.
\end{eqnarray}
\end{Theorem}

This criterion was applied to tensor products of Euler and Wiener integrated processes (see \cite{LifPapWoz2}) and to the tensor product approximation problems with weighted Korobov kernels (see \cite{LifPapWoz1} and \cite{Xu}). However, the practice shows that the criterion is not convenient enough for applicaions to concrete examples of $(X_d)_{d\in\N}$, because verification of  \eqref{th_oldQPTcrit_cond} usually requires a lot of additional calculations. The next theorem states that \eqref{th_oldQPTcrit_cond} can be splited into two simpler conditions of additive type, where one of them is close to \eqref{th_oldQPTcrit_necesscond}.

\begin{Theorem}\label{th_newQPTcrit}
The sequence of approximation problems for $X_d=\otimes_{j=1}^d X_{1,j}$, $d\in\N$, is quasi-polynomially tractable iff the following both conditions are verified
\begin{eqnarray}
&&\sup_{d\in\N} \dfrac{1}{\ln_+ d}\sum_{j=1}^{d}\sum_{k=2}^{\infty} \bigl(1+ |\ln\nlambda^{X_{1,j}}_k|\bigr)\,\nlambda^{X_{1,j}}_k<\infty,\label{th_newQPTcrit_ln}\\
&&\sup_{d\in\N} \sum_{j=1}^{d} \sum_{k=2}^{\infty} \bigl(\nlambda^{X_{1,j}}_k\bigr)^{1-\tfrac{\gamma}{\ln_+ d}}\,\id\Bigl(\nlambda^{X_{1,j}}_k<e^{-\tau\ln_+ d}\Bigr)<\infty\label{th_newQPTcrit_exp}
\end{eqnarray}
for some $\gamma\in(0,1)$ and for some $($or equivalently each\,$)$ $\tau>0$.
\end{Theorem}
\textbf{Proof of Theorem \ref{th_newQPTcrit}.}\quad \textit{Sufficiency}.  Suppose that we have \eqref{th_newQPTcrit_ln} and \eqref{th_newQPTcrit_exp} for some $\gamma\in(0,1)$ and $\tau>0$. Without loss of generality we assume that $\tau\geqslant 1$. Let us consider the products from the condition \eqref{th_oldQPTcrit_cond}. Using \eqref{cond_sumnlambdaX1jk}, we infer the following representation
\begin{eqnarray}\label{conc_prodQPT_repres}
\prod_{j=1}^{d} \sum_{k=1}^{\infty} \bigl(\nlambda^{X_{1,j}}_k\bigr)^{1-\tfrac{\gamma}{\ln_+ d}}= \prod_{j=1}^{d} \biggl(1+\sum_{k=1}^{\infty}\Bigl( \bigl(\nlambda^{X_{1,j}}_k\bigr)^{1-\tfrac{\gamma}{\ln_+ d}}-\nlambda^{X_{1,j}}_k\Bigr)\biggr).
\end{eqnarray}
Next, applying elementary inequality $1+x<e^x$, $x\geqslant 0$, we obtain
\begin{eqnarray}
\prod_{j=1}^{d} \sum_{k=1}^{\infty} \bigl(\nlambda^{X_{1,j}}_k\bigr)^{1-\tfrac{\gamma}{\ln_+ d}}
&\leqslant& \exp\biggl\{\sum_{j=1}^d \sum_{k=1}^{\infty}\Bigl( \bigl(\nlambda^{X_{1,j}}_k\bigr)^{1-\tfrac{\gamma}{\ln_+ d}}-\nlambda^{X_{1,j}}_k\Bigr)\biggr\}\nonumber\\
&\leqslant& \exp\biggl\{\sum_{j=1}^d S^{X_{1,j}}_{d,\gamma}(\tau)+\sum_{j=1}^d R^{X_{1,j}}_{d,\gamma}(\tau)\biggr\},\label{conc_prodQPT_ineq}
\end{eqnarray}
where we set
\begin{eqnarray}
S^{X_{1,j}}_{d,\gamma}(\tau)&\colonequals&\sum_{k=1}^\infty \Bigl(\bigl(\nlambda^{X_{1,j}}_k\bigr)^{1-\tfrac{\gamma}{\ln_+ d}}-\nlambda^{X_{1,j}}_k\Bigr)\id\Bigl(\nlambda^{X_{1,j}}_k\geqslant e^{-\tau\ln_+d}\Bigr),\nonumber\\
R^{X_{1,j}}_{d,\gamma}(\tau)&\colonequals&\sum_{k=1}^\infty \bigl(\nlambda^{X_{1,j}}_k\bigr)^{1-\tfrac{\gamma}{\ln_+ d}}\id\Bigl(\nlambda^{X_{1,j}}_k< e^{-\tau\ln_+d}\Bigr).\label{def_RX1jgamtau}
\end{eqnarray}
We first consider the sums $S^{X_{1,j}}_{d,\gamma}(\tau)$, $d\in\N$. Let us find the constant $C_{\gamma,\tau}$ such that $e^{\gamma  x}\leqslant1+C_{\gamma,\tau} x $ for any $x\in[0,\tau]$. Using this inequality we estimate
\begin{eqnarray*}
S^{X_{1,j}}_{d,\gamma}(\tau)&=&\sum_{k=1}^\infty \biggl( \exp\Bigl\{\gamma \tfrac{|\ln \nlambda^{X_{1,j}}_k |}{\ln_+ d}\Bigr\}-1\biggr)\nlambda^{X_{1,j}}_k\id\Bigl(\nlambda^{X_{1,j}}_k\geqslant e^{-\tau\ln_+d}\Bigr)\\
&\leqslant&\dfrac{C_{\gamma,\tau}}{\ln_+ d}\sum_{k=1}^{\infty} \bigl|\ln\nlambda^{X_{1,j}}_k\bigr|\,\nlambda^{X_{1,j}}_k .
\end{eqnarray*}
On account of the inequality $\ln(1+x)\leqslant x$, $x\geqslant 0$, and \eqref{cond_sumnlambdaX1jk}, observe that
\begin{eqnarray*}
|\ln\nlambda^{X_{1,j}}_1|\,\nlambda^{X_{1,j}}_1=\ln\biggl(1+\tfrac{1-\nlambda^{X_{1,j}}_1}{\nlambda^{X_{1,j}}_1}\biggl)\,\nlambda^{X_{1,j}}_1\leqslant 1-\nlambda^{X_{1,j}}_1=\sum_{k=2}^{\infty} \nlambda^{X_{1,j}}_k.
\end{eqnarray*}
Therefore
\begin{eqnarray*}
S^{X_{1,j}}_{d,\gamma}(\tau)\leqslant  \dfrac{C_{\gamma,\tau}}{\ln_+ d}\sum_{k=2}^{\infty} \bigl(1+|\ln\nlambda^{X_{1,j}}_k|\bigr)\,\nlambda^{X_{1,j}}_k.
\end{eqnarray*}
Thus from this and  \eqref{th_newQPTcrit_ln} we conclude that 
\begin{eqnarray}\label{conc_SX1jdgamtau}
\sup_{d\in\N} \sum_{j=1}^d S^{X_{1,j}}_{d,\gamma}(\tau)<\infty.
\end{eqnarray}

We next consider the sums $R^{X_{1,j}}_{d,\gamma}(\tau)$, $d\in\N$:
\begin{eqnarray*}
R^{X_{1,j}}_{d,\gamma}(\tau)\leqslant d^{-\tau}e^{\gamma\tau}+ \sum_{k=2}^\infty \bigl(\nlambda^{X_{1,j}}_k\bigr)^{1-\tfrac{\gamma}{\ln_+ d}}\id\Bigl(\nlambda^{X_{1,j}}_k< e^{-\tau\ln_+d}\Bigr).
\end{eqnarray*}
Since $\tau\geqslant1$, it follows that
\begin{eqnarray*}
\sum_{j=1}^{d}R^{X_{1,j}}_{d,\gamma}(\tau)\leqslant e^{\gamma\tau}+ \sum_{j=1}^{d} \sum_{k=2}^\infty \bigl(\nlambda^{X_{1,j}}_k\bigr)^{1-\tfrac{\gamma}{\ln_+ d}}\id\Bigl(\nlambda^{X_{1,j}}_k< e^{-\tau\ln_+d}\Bigr).
\end{eqnarray*}
According to \eqref{th_newQPTcrit_exp} we obtain
\begin{eqnarray}\label{conc_RX1jdgamtau}
\sup_{d\in\N} \sum_{j=1}^d R^{X_{1,j}}_{d,\gamma}(\tau)<\infty.
\end{eqnarray}
Thus we conclude from \eqref{conc_prodQPT_ineq}, \eqref{conc_SX1jdgamtau}, and \eqref{conc_RX1jdgamtau} that the condition \eqref{th_oldQPTcrit_cond} of Theorem \ref{th_oldQPTcrit} holds for given $\gamma$. Hence we have the quasi-polynomial tractability.

\textit{Necessity}. Suppose that the sequence of approximation problems for $X_d=\otimes_{j=1}^d X_{1,j}$, $d\in\N$, is quasi-polynomially tractable.  

We first show that \eqref{th_newQPTcrit_ln} is satisfied. On the one hand, the quantity $n^{X_d}(\e)$ satisfies \eqref{def_QPT} for some $C>0$ and $s\geqslant0$. On the other hand, from \eqref{def_nXde}, \eqref{conc_LambdaXd}, and $\lambda^{X_d}_1=\prod_{j=1}^d \lambda^{X_{1,j}}_1$ we have the inequality
\begin{eqnarray*}
n^{X_d}(\e)\geqslant (1-\e^2) \dfrac{\Lambda^{X_d}}{\lambda^{X_d}_1} = (1-\e^2) \prod_{j=1}^d \dfrac{1}{\nlambda^{X_{1,j}}_1}\quad\text{for all}\quad d\in\N,\,\, \e\in(0,1). 
\end{eqnarray*}
Consequently, 
\begin{eqnarray*}
\ln(1-\e^2)+\sum_{j=1}^d |\ln \nlambda^{X_{1,j}}_1|
\leqslant \ln C + s(1+|\ln \e|)(1+\ln d).
\end{eqnarray*}
Hence for all $d\in\N$
\begin{eqnarray*}
\sup_{d\in\N}\dfrac{1}{\ln_+ d}\sum_{j=1}^d |\ln \nlambda^{X_{1,j}}_1|<\infty.
\end{eqnarray*}
Applying the elementary inequality $\ln x \leqslant x-1$, $x>0$, and \eqref{cond_sumnlambdaX1jk}, observe that  
\begin{eqnarray*}
|\ln \nlambda^{X_{1,j}}_1|\geqslant 1-\nlambda^{X_{1,j}}_1=\sum_{k=2}^\infty \nlambda^{X_{1,j}}_k,\quad j\in\N.
\end{eqnarray*}
Therefore
\begin{eqnarray*}
\sup_{d\in\N} \dfrac{1}{\ln_+ d}\sum_{j=1}^d\sum_{k=2}^\infty \nlambda^{X_{1,j}}_k<\infty.
\end{eqnarray*}
On account of necessary condition \eqref{th_oldQPTcrit_necesscond} from Theorem \ref{th_oldQPTcrit} we obtain \eqref{th_newQPTcrit_ln}.

We next prove \eqref{th_newQPTcrit_exp} for some $\gamma\in(0,1)$ and for any $\tau>0$. By Theorem \ref{th_oldQPTcrit}, we have \eqref{th_oldQPTcrit_cond} for some $\gamma=\gamma_*\in(0,1)$. From the representation \eqref{conc_prodQPT_repres} we conclude
\begin{eqnarray*}
\prod_{j=1}^{d} \sum_{k=1}^{\infty} \bigl(\nlambda^{X_{1,j}}_k\bigr)^{1-\tfrac{\gamma_*}{\ln_+ d}}
\geqslant 1+\sum_{j=1}^d \sum_{k=1}^{\infty}\Bigl( \bigl(\nlambda^{X_{1,j}}_k\bigr)^{1-\tfrac{\gamma_*}{\ln_+ d}}-\nlambda^{X_{1,j}}_k\Bigr).
\end{eqnarray*}
Thus we have
\begin{eqnarray}\label{conc_prodQPT_sumexp}
\sup_{d\in\N}\sum_{j=1}^d \sum_{k=1}^{\infty}\Bigl( \bigl(\nlambda^{X_{1,j}}_k\bigr)^{1-\tfrac{\gamma_*}{\ln_+ d}}-\nlambda^{X_{1,j}}_k\Bigr)<\infty.
\end{eqnarray}
Choose any $\gamma\in(0,\gamma_*)$ and any $\tau>0$. Consider the sum $R^{X_{1,j}}_{d,\gamma}(\tau)$, which is defined by \eqref{def_RX1jgamtau}. It admits the following integral representation
\begin{eqnarray*}
R^{X_{1,j}}_{d,\gamma}(\tau)=-\int\limits_\tau^\infty e^{\gamma t}\dd R^{X_{1,j}}_{d,0}(t).
\end{eqnarray*}
Integrating by parts yields
\begin{eqnarray}\label{conc_RX1jdgamtau2}
R^{X_{1,j}}_{d,\gamma}(\tau)=e^{\gamma \tau} R^{X_{1,j}}_{d,0}(\tau)-\lim_{t\to\infty} e^{\gamma t} R^{X_{1,j}}_{d,0}(t) +\gamma\int\limits_\tau^\infty e^{\gamma x} R^{X_{1,j}}_{d,0}(t) \dd x.
\end{eqnarray}
It is easy to prove that
\begin{eqnarray*}
R^{X_{1,j}}_{d,0}(t)=\sum_{k=1}^{\infty} \nlambda^{X_{1,j}}_k\id\bigl(\nlambda^{X_{1,j}}_k< e^{-t\ln_+ d}\bigr)\leqslant \dfrac{1}{e^{\gamma_* t}-1}\sum_{k=1}^\infty \Bigl(\bigl(\nlambda^{X_{1,j}}_k\bigr)^{1-\tfrac{\gamma_*}{\ln_+ d}}-\nlambda^{X_{1,j}}_k\Bigr).
\end{eqnarray*}
From this inequality we conclude that the limit in the previous expression in \eqref{conc_RX1jdgamtau2} exists and equals zero. Applying this inequality to other terms of \eqref{conc_RX1jdgamtau2} we get
\begin{eqnarray*}
R^{X_{1,j}}_{d,\gamma}(\tau)\leqslant M_{\gamma, \tau}\sum_{k=1}^\infty \Bigl(\bigl(\nlambda^{X_{1,j}}_k\bigr)^{1-\tfrac{\gamma_*}{\ln_+ d}}-\nlambda^{X_{1,j}}_k\Bigr),
\end{eqnarray*}
where we set
\begin{eqnarray*}
M_{\gamma,\tau}\colonequals\dfrac{e^{\gamma\tau}}{e^{\gamma_*\tau}-1} +\gamma\int\limits_\tau^\infty  \dfrac{e^{\gamma t}}{e^{\gamma_* t}-1}\dd t<\infty.
\end{eqnarray*}
From this we conclude that
\begin{eqnarray*}
\sum_{j=1}^d R^{X_{1,j}}_{d,\gamma}(\tau)\leqslant M_{\gamma, \tau}\sum_{j=1}^d\sum_{k=1}^\infty \Bigl(\bigl(\nlambda^{X_{1,j}}_k\bigr)^{1-\tfrac{\gamma_*}{\ln_+ d}}-\nlambda^{X_{1,j}}_k\Bigr).
\end{eqnarray*}
In view of \eqref{conc_prodQPT_sumexp} we have \eqref{th_newQPTcrit_exp}.\quad$\Box$\\

We comment on the conditions of Theorem \ref{th_newQPTcrit}. Typically, for concrete examples of $(X_d)_{d\in\N}$ only \eqref{th_newQPTcrit_ln} is important for quasi-polynomial tractability, because the condition \eqref{th_newQPTcrit_exp} usually holds under the natural assumptions on the sequence $(X_d)_{d\in\N}$. As we will see below, the verification of \eqref{th_newQPTcrit_exp} is rather simple.

\section{Applications}
\subsection{Korobov kernels}\newparagraph
Let $B_{g,r}(t)$, $t\in[0,1]$, be a zero-mean random process with the following covariance function
\begin{eqnarray*}
\CorrFunc^{B_{g,r}}(t,s)\colonequals 1+2g\sum_{k=1}^{\infty} k^{-2r}\cos(2\pi k(t-s)), \quad t,s\in[0,1],
\end{eqnarray*} 
which is called \textit{Korobov kernel}. Here $g\in(0,1]$ and $r>1/2$.

We consider $B_{g,r}(t)$, $t\in[0,1]$, as a random element $B_{g,r}$ of the space $L_2([0,1])$. The covariance operator $K^{B_{r,g}}$ of $B_{r,g}$ is an integration operator with kernel $\CorrFunc^{B_{g,r}}$. The eigenvalues of $K^{B_{r,g}}$ are exactly known (see \cite{NovWoz1}):
\begin{eqnarray}\label{def_Koreigen}
\lambda^{B_{g,r}}_1=1,\quad  \lambda^{B_{g,r}}_{2k}=\lambda^{B_{g,r}}_{2k+1}=\dfrac{g}{k^{2r}},\quad k\in\N.
\end{eqnarray}
Note that the trace of $K^{B_{r,g}}$ is
\begin{eqnarray*}
\Lambda^{B_{g,r}}=1+2g \zeta(2r),
\end{eqnarray*}
where $\zeta(p)=\sum_{k=1}^\infty k^{-p}$, $p>1$, is the Riemann zeta-function.

Suppose that we have a sequence of processes $B_{g_j,r_j}(t)$, $t\in[0,1]$,  with covariance functions $\CorrFunc^{B_{g_j,r_j}}$, $j\in\N$, respectively. Let $\KorElem_d(t)$, $t\in[0,1]^d$, $d\in\N$, be the sequence of zero-mean random fields with the following covariance functions
\begin{eqnarray*}
\CorrFunc^{\KorElem_d}(t,s)=\prod_{j=1}^{d}\CorrFunc^{B_{g_j,r_j}}(t_j,s_j), \quad t,s\in[0,1]^d,\quad d\in\N.
\end{eqnarray*}
Every field $\KorElem_d(t)$, $t\in\R^d$, can be considered as a random element $\KorElem_d$ of the space $L_2([0,1]^d)$. Every $\KorElem_d$ has a covariance operator of the tensor product form $K^{\KorElem_d}=\otimes_{j=1}^d K^{B_{g_j,r_j}}$, $d\in\N$, i.e. by definition from Section 2, $\KorElem_d=\otimes_{j=1}^d B_{g_j,r_j}$, $d\in\N$.

In \cite{LifPapWoz1} M. A. Lifshits, A. Papageorgiou and H. Wo\'zniakowski were the first to investigate approximation problems for $\KorElem_d$, $d\in\N$, in the average case setting. Under the assumptions
\begin{eqnarray}\label{cond_rjgj}
1\geqslant g_1\geqslant g_2\geqslant\ldots\geqslant g_j\geqslant\ldots >0,\qquad 1/2<r_1\leqslant r_2\leqslant \ldots\leqslant r_j\leqslant\ldots,
\end{eqnarray}
they proved that quasi-polynomial tractability holds whenever \eqref{th_KorQPT_cond} is satisfied (see below) and $\liminf_{j\to\infty} (r_j/\ln j)>0$. In the recent paper \cite{Xu} G. Xu shows that the latter condition can be omitted under \eqref{cond_rjgj}. Moreover, the next theorem asserts that there is no need to assume monotonicity for $(g_j)_{j\in\N}$ and $(r_j)_{j\in\N}$.

\begin{Theorem}\label{th_KorQPT}
Let $(r_j)_{j\in\N}$ be a sequence such that $\inf_{j\in\N} r_j>1/2$. Let $(g_j)_{j\in\N}$ be a positive sequence such that $\sup_{j\in\N} g_j\leqslant 1$. The sequence of approximation problems for $\KorElem_d=\otimes_{j=1}^d B_{g_j,r_j}$, $d\in\N$, is quasi-polynomially tractable iff 
\begin{eqnarray}\label{th_KorQPT_cond}
\sup_{d\in\N} \dfrac{1}{\ln_+ d}\sum_{j=1}^d (1+|\ln g_j|)\, g_j<\infty.
\end{eqnarray}
\end{Theorem}
\textbf{Proof of Theorem \ref{th_KorQPT}.}\quad  Define $r_0\colonequals \inf_{j\in\N} r_j>1/2$ and $g_0\colonequals \sup_{j\in\N} g_j\leqslant 1$. For every $j\in\N$ we consider the following sum 
\begin{eqnarray*}
L_{g_j,r_j}\colonequals \sum_{k=2}^{\infty} \bigl(1+ |\ln\nlambda^{B_{g_j,r_j}}_k|\bigr)\,\nlambda^{B_{g_j,r_j}}_k
=2\sum_{k=1}^{\infty} \dfrac{\bigl(1+\ln(1+2g_j \zeta(2r_j))+\ln (k^{2r_j})+|\ln g_j|\bigr)\,g_j}{k^{2r_j}(1+2g_j \zeta(2r_j)) }
\end{eqnarray*}
from the condition \eqref{th_newQPTcrit_ln} of Theorem \ref{th_newQPTcrit}, where we set $X_{1,j}=B_{g_j,r_j}$, $j\in\N$. We first provide the lower estimate for every $L_{g_j,r_j}$:
\begin{eqnarray*}
L_{g_j,r_j}\geqslant 2C_1(1+|\ln g_j|)\, g_j, \quad j\in\N,
\end{eqnarray*}
where $C_1\colonequals (1+2g_0 \zeta(2r_0))^{-1}$. Next, we obtain the upper estimate for $L_{g_j,r_j}$:
\begin{eqnarray*}
L_{g_j,r_j}&\leqslant& 2\sum_{k=1}^{\infty} \bigl(1+\ln(1+2g_0 \zeta(2r_0))+\ln (k^{2r_j})+|\ln g_j|\bigr)\,\dfrac{g_j}{k^{2r_j}}\\
&\leqslant& 2(C_2+C_3+C_4)(1+|\ln g_j|)\, g_j,\quad j\in\N,
\end{eqnarray*}
where $C_2\colonequals \zeta(2r_0)$, $C_3\colonequals 2 g_0 \zeta(2r_0)^2$, $C_4\colonequals  \sup_{j\in\N}\sum_{k=1}^\infty \tfrac{\ln (k^{2r_j})}{k^{2r_j}}$. Thus we have
\begin{eqnarray*}
\dfrac{1}{\ln_+ d}\sum_{j=1}^{d} L_{g_j,r_j}\asymp \dfrac{1}{\ln_+ d}\sum_{j=1}^{d} (1+|\ln g_j|)\, g_j, \quad d\in\N,
\end{eqnarray*}
i.e. the condition \eqref{th_newQPTcrit_ln} of Theorem \ref{th_newQPTcrit} for $B_{g_j,r_j}$ is equivalent to \eqref{th_KorQPT_cond}.

Next, we verify that the condition \eqref{th_newQPTcrit_exp} of Theorem \ref{th_newQPTcrit} for $B_{g_j,r_j}$, $j\in\N$, holds for some $\gamma\in(0,1)$ and $\tau>0$. Fix $\gamma\in(0,1)$ such that $2r_0(1-\gamma)>1$ and consider the quantity
\begin{eqnarray*}
R^{B_{g_j,r_j}}_{d,\gamma}(\tau)\colonequals\sum_{k=2}^{\infty} \bigl(\nlambda^{B_{g_j,r_j}}_k\bigr)^{1-\tfrac{\gamma}{\ln_+ d}}\,\id\Bigl(\nlambda^{B_{g_j,r_j}}_k<e^{-\tau\ln_+ d}\Bigr)=2\sum_{k=k_{d,j}(\tau)}^{\infty} \biggl(\dfrac{g_j k^{-2r_j}}{1+2g_j \zeta(2r_j)}\biggr)^{1-\tfrac{\gamma}{\ln_+d}},
\end{eqnarray*}
where we set
\begin{eqnarray}\label{def_kdjtau}
k_{d,j}(\tau)\colonequals \min\biggl\{k\in\N: \dfrac{g_j k^{-2r_j}}{1+2g_j \zeta(2r_j)}<e^{-\tau\ln_+d} \biggr\}.
\end{eqnarray}
From this we infer the following inequality
\begin{eqnarray}\label{conc_kdjtau}
k_{d,j}(\tau)-1\leqslant \biggl(\dfrac{g_j e^{\tau\ln_+d}}{1+2g_j \zeta(2r_j)}\biggr)^{\tfrac{1}{2r_j}}\leqslant \bigl(g_0 e^{\tau\ln_+d}\bigr)^{\tfrac{1}{2r_0}}.
\end{eqnarray}
Using $\sum_{k=n}^{\infty} f(k)\leqslant f(n)+\int_n^\infty f(t)\dd t$ for monotonic non-increasing $f$, we estimate
\begin{eqnarray*}
R^{B_{g_j,r_j}}_{d,\gamma}(\tau)&\leqslant& 2\biggl(\dfrac{g_j k_{d,j}(\tau)^{-2r_j}}{1+2g_j \zeta(2r_j)}\biggr)^{1-\tfrac{\gamma}{\ln_+d}}+2\int\limits_{k_{d,j}(\tau)}^\infty \biggl(\dfrac{g_j t^{-2r_j}}{1+2g_j \zeta(2r_j)}\biggr)^{1-\tfrac{\gamma}{\ln_+d}}\dd t\\
&=& 2\biggl(1+\dfrac{k_{d,j}(\tau)}{2r_j\bigl(1-\tfrac{\gamma}{\ln_+d}\bigr) -1}\biggr)\biggl(\dfrac{g_j k_{d,j}(\tau)^{-2r_j}}{1+2g_j \zeta(2r_j)}\biggr)^{1-\tfrac{\gamma}{\ln_+d}}.
\end{eqnarray*}
According to \eqref{def_kdjtau} and \eqref{conc_kdjtau} we have
\begin{eqnarray*}
R^{B_{g_j,r_j}}_{d,\gamma}(\tau)\leqslant 2\Biggl(1+\dfrac{1+g_0^{\tfrac{1}{2r_0}}\cdot e^{\tfrac{\tau\ln_+d}{2r_0}}}{2r_j(1-\gamma) -1}\Biggr)\bigl(e^{-\tau\ln_+d}\bigr)^{1-\tfrac{\gamma}{\ln_+d}}\leqslant C_5 \exp\Bigl\{ -\tau\bigl(1-\tfrac{1}{2r_0}\bigr)\ln_+ d\Bigr\},
\end{eqnarray*}
where $C_5\colonequals 2\cdot\tfrac{2r_0(1-\gamma)+g_0^{1/2r_0}}{2r_0(1-\gamma)-1}\cdot e^{\gamma\tau}$.

Next, choose any $\tau$ such that $\tau\bigl(1-\tfrac{1}{2r_0}\bigr)\geqslant 1$. Then
\begin{eqnarray*}
\sum_{j=1}^d R^{B_{g_j,r_j}}_{d,\gamma}(\tau)\leqslant d\cdot C_5 \exp\Bigl\{ -\tau\bigl(1-\tfrac{1}{2r_0}\bigr)\ln_+ d\Bigr\}\leqslant C_5,\quad d\in\N.
\end{eqnarray*}
Hence \eqref{th_newQPTcrit_exp} holds for $B_{g_j,r_j}$, $j\in\N$. Thus \eqref{th_KorQPT_cond} is necessary and sufficient condition for quasi-polynomial tractability \quad $\Box$.

\subsection{Squared exponential kernels}\newparagraph
Let $G_\sigma(t)$, $t\in\R$, be a zero-mean random process with the following covariance function
\begin{eqnarray*}
\CorrFunc^{G_\sigma}(t,s)\colonequals e^{-\tfrac{(t-s)^2}{2\sigma^2}}, \quad t,s\in\R,
\end{eqnarray*} 
where $\sigma>0$ is a characteristic length-scale. The function $\CorrFunc^{G_\sigma}$ is rather popular kernel function used in machine learning (see \cite{RasWill}).
We consider $G_\sigma(t)$, $t\in\R$, as a random element $G_\sigma$ of the space $L_2(\R,\mu)$, where $\mu$ is a standard Gaussian distribution on $\R$. The covariance operator $K^{G_\sigma}$ of $G_\sigma$ acts as follows
\begin{eqnarray*}
K^{G_\sigma}f(t)=\int\limits_{\R} \CorrFunc^{G_\sigma}(t,s) f(s)\mu(\!\dd s)=\int\limits_{\R} e^{-\tfrac{(t-s)^2}{2\sigma^2}} f(s)\, \dfrac{1}{\sqrt{2\pi}}\, e^{-\tfrac{s^2}{2}} \dd s,\quad t\in\R.
\end{eqnarray*}
Eigenvalues of $K^{G_\sigma}$ are known (see  \cite{NovWoz3} and \cite{RasWill}):
\begin{eqnarray}\label{def_SEeigen}
\lambda^{G_\sigma}_k=(1-\omega_\sigma)\,\omega_\sigma^{k-1}, \quad k\in\N,
\end{eqnarray}
where $\omega_\sigma\colonequals (1+\sigma^2 I_\sigma)^{-1}$, $I_\sigma\colonequals \tfrac{1}{2}+\tfrac{1}{2} \sqrt{1+\tfrac{4}{\sigma^2}}$. It is easily seen that $\Lambda^{G_\sigma}=\sum_{k\in\N} \lambda^{G_\sigma}_k=1$, i.e. $\lambda^{G_\sigma}_k=\nlambda^{G_\sigma}_k$, $k\in\N$.

Suppose that we have the sequence of processes $G_{\sigma_j}(t)$, $t\in\R$,  with covariance functions $\CorrFunc^{G_{\sigma_j}}$, $j\in\N$, respectively. Consider the sequence of zero-mean random fields $\SEElem_d(t)$, $t\in\R^d$, $d\in\N$, with the following covariance functions
\begin{eqnarray*}
\CorrFunc^{\SEElem_d}(t,s)=\prod_{j=1}^{d}\CorrFunc^{G_{\sigma_j}}(t_j,s_j), \quad t,s\in\R^d,\quad d\in\N.
\end{eqnarray*}
Every field $\SEElem_d(t)$, $t\in\R^d$, is a random element $\SEElem_d$ of the space $L_2(\R^d,\mu^d)$, where $\mu^d$ is a standard Gaussian measure on $\R^d$. Thus we have $\SEElem_d=\otimes_{j=1}^d G_{\sigma_j}$, $d\in\N$. We find the criterion of quasi-polynomial tractability of approximation problems for these elements (worst case setting results can be found in \cite{FasshHick} and \cite{NovWoz3}).

\begin{Theorem}\label{th_SEQPT}
Let $(\sigma_j)_{j\in\N}$ be a sequence such that $\inf_{j\in\N} \sigma_j>0$.  The sequence of approximation problems for $\SEElem_d=\otimes_{j=1}^d G_{\sigma_j}$, $d\in\N$, is quasi-polynomially tractable iff
\begin{eqnarray}\label{th_SEQPT_cond}
\sup_{d\in\N} \dfrac{1}{\ln_+ d}\sum_{j=1}^d \dfrac{1+\ln(1+ \sigma_j^2)}{\sigma_j^2}<\infty.
\end{eqnarray}
\end{Theorem}
\textbf{Proof of Theorem \ref{th_SEQPT}.}\quad Let us consider the sums
\begin{eqnarray*}
L_{\sigma_j}\colonequals \sum_{k=2}^{\infty} \bigl(1+ |\ln\nlambda^{G_{\sigma_j}}_k|\bigr)\,\nlambda^{G_{\sigma_j}}_k,\quad j\in\N,
\end{eqnarray*}
from the condition \eqref{th_newQPTcrit_ln} of Theorem \ref{th_newQPTcrit}, where we set $X_{1,j}=G_{\sigma_j}$. Using \eqref{def_SEeigen} we find
\begin{eqnarray*}
L_{\sigma_j}&=& \sum_{k=2}^{\infty} \bigl( 1-\ln(1-\omega_{\sigma_j})-(k-1)\ln \omega_{\sigma_j}\bigr) (1-\omega_{\sigma_j})\, \omega_{\sigma_j}^{k-1}\\
&=&\bigl( 1-\ln(1-\omega_{\sigma_j})\bigr)(1-\omega_{\sigma_j}) \sum_{k=2}^{\infty} \omega_{\sigma_j}^{k-1}\\
&&{}-\ln (\omega_{\sigma_j}) (1-\omega_{\sigma_j})\, \omega_{\sigma_j}\sum_{k=2}^{\infty} (k-1)\, \omega_{\sigma_j}^{k-2}\\
&=& \bigl(1-\ln(1-\omega_{\sigma_j})\bigr)\,\omega_{\sigma_j}-\dfrac{\ln (\omega_{\sigma_j})\,\omega_{\sigma_j}}{1-\omega_{\sigma_j}}.
\end{eqnarray*}
Substituting $\omega_{\sigma_j}= (1+\sigma_j^2 I_{\sigma_j})^{-1}$ in the last representation for $L_{\sigma_j}$, we infer
\begin{eqnarray}\label{conc_Lsigmaj}
L_{\sigma_j}=\dfrac{1+\ln\Bigl(1+\tfrac{1}{\sigma_j^2I_{\sigma_j}}\Bigr)}{1+\sigma_j^2 I_{\sigma_j}}+ \dfrac{\ln\bigl(1+\sigma_j^2I_{\sigma_j}\bigr)}{\sigma_j^2 I_{\sigma_j}}.
\end{eqnarray}
Define $\sigma_0:=\inf_{j\in\N} \sigma_j>0$. For any $j\in\N$ we have $1<I_{\sigma_j}\leqslant I_{\sigma_0}$. Using the inequality $\ln(1+x)\leqslant x$, $x\geqslant0$, we see that
\begin{eqnarray*}
L_{\sigma_j}
\leqslant\dfrac{1}{\sigma_j^2 I_{\sigma_j}}+ \dfrac{\ln\bigl(1+\sigma_j^2 I_{\sigma_j}\bigr)}{\sigma_j^2 I_{\sigma_j}}
\leqslant\dfrac{1+\ln\bigl(1+\sigma_j^2I_{\sigma_0}\bigr)}{\sigma_j^2}.
\end{eqnarray*}
Consequently, we have the following upper estimate
\begin{eqnarray*}
L_{\sigma_j}\leqslant C_1 \dfrac{1+\ln\bigl(1+\sigma_j^2\bigr)}{\sigma_j^2},\quad j\in\N, 
\end{eqnarray*}
where $C_1\colonequals 1+\ln I_{\sigma_0}>0$. Next, we see that by \eqref{conc_Lsigmaj}
\begin{eqnarray*}
L_{\sigma_j}\geqslant\dfrac{1+\ln\bigl(1+\sigma_j^2I_{\sigma_j}\bigr)}{1+\sigma_j^2 I_{\sigma_j}}
\geqslant\dfrac{1+\ln\bigl(1+\sigma_j^2\bigr)}{1+\sigma_j^2 I_{\sigma_0}}.
\end{eqnarray*}
Hence we obtain the following lower estimate
\begin{eqnarray*}
L_{\sigma_j}\geqslant C_2 \dfrac{1+\ln\bigl(1+\sigma_j^2\bigr)}{\sigma_j^2},\quad j\in\N, 
\end{eqnarray*}
where $C_2\colonequals \sigma_0^2/(1+\sigma_0^2I_{\sigma_0})>0$. Therefore 
\begin{eqnarray*}
\dfrac{1}{\ln_+ d}\sum_{j=1}^{d} L_{\sigma_j}\asymp \dfrac{1}{\ln_+ d}\sum_{j=1}^{d} \dfrac{1+\ln\bigl(1+\sigma_j^2\bigr)}{\sigma_j^2}, \quad d\in\N.
\end{eqnarray*}
Thus the condition \eqref{th_newQPTcrit_ln} of Theorem \ref{th_newQPTcrit} is equivalent to \eqref{th_SEQPT_cond}.

Next, we verify that for $(G_{\sigma_j})_{j\in\N}$ the condition \eqref{th_newQPTcrit_exp} of Theorem \ref{th_newQPTcrit} is always satisfied under the assumption $\sigma_0>0$. Fix any $\gamma\in(0,1)$ and consider the quantity
\begin{eqnarray*}
R^{G_{\sigma_j}}_{d,\gamma}\colonequals\sum_{k=2}^{\infty} \bigl(\nlambda^{G_{\sigma_j}}_k\bigr)^{1-\tfrac{\gamma}{\ln_+ d}}\,\id\Bigl(\nlambda^{G_{\sigma_j}}_k<e^{-\ln_+ d}\Bigr),\quad j\in\N.
\end{eqnarray*}
Let us introduce the threshold index
\begin{eqnarray*}
k_{d,j}\colonequals \min\{k\in\N: (1-\omega_{\sigma_j})\,\omega_{\sigma_j}^{k-1} <e^{-\ln_+ d}, k\geqslant 2\}.
\end{eqnarray*}
According to \eqref{def_SEeigen} we infer
\begin{eqnarray*}
R^{G_{\sigma_j}}_{d,\gamma}(\tau)
=\sum_{k=k_{d,j}}^{\infty} \bigl((1-\omega_{\sigma_j})\,\omega_{\sigma_j}^{k-1}\bigr)^{1-\tfrac{\gamma}{\ln_+ d}}
=\dfrac{\bigl((1-\omega_{\sigma_j})\,\omega_{\sigma_j}^{k_{d,j}-1}\bigr)^{1-\tfrac{\gamma}{\ln_+ d}}}{1- \omega_{\sigma_j}^{1-\tfrac{\gamma}{\ln_+ d}}}.
\end{eqnarray*}
By definition $k_{d,j}$, we see that
\begin{eqnarray*}
R^{G_{\sigma_j}}_{d,\gamma}(\tau)\leqslant \dfrac{\bigl((1-\omega_{\sigma_j})\,\omega_{\sigma_j}^{k_{j,d}-1}\bigr)^{1-\tfrac{\gamma}{\ln_+ d}}}{1-\omega_{\sigma_0}^{1-\gamma}} \leqslant \dfrac{\bigl(e^{-\ln_+d}\bigr)^{1-\tfrac{\gamma}{\ln_+ d}}}{1-\omega_{\sigma_0}^{1-\gamma}}\leqslant \dfrac{e^{\gamma-\ln_+d}}{1-\omega_{\sigma_0}^{1-\gamma}}.
\end{eqnarray*}
Then $\sum_{j=1}^{d} R^{G_{\sigma_j}}_{d,\gamma}\leqslant e^{\gamma}/(1-\omega_{\sigma_0}^{1-\gamma})$ for any $d\in\N$. Thus \eqref{th_newQPTcrit_exp} holds as required.\quad $\Box$.

\section*{Acknowlegments}\newparagraph
The author is grateful to Professor M. A. Lifshits for suggesting the problem and for several helpful comments.

\textit{Keywords and phrases}: linear tensor product approximation problems, average case approximation complexity, quasi-polynomial tractability, random fields.

\textsc{Department of Mathematics and Mechanics,  St. Petersburg State University, Universitetsky pr. 28, 198504 St. Petersburg, Russia.}\\
\textit{E-mail address}: \texttt{alexeykhartov@gmail.com} 
\end{document}